\documentclass[letterpaper, 10 pt, conference]{ieeeconf}  

\IEEEoverridecommandlockouts 
\usepackage{graphics} 
\usepackage{amsmath} 
\usepackage{amssymb}  
\usepackage{mathtools}
\usepackage{hyperref}
\let\labelindent\relax
\usepackage{enumitem}

\usepackage{amsmath}

\newtheorem{proposition}{Proposition}
\usepackage[font=footnotesize]{caption}
\usepackage{cite}
\usepackage{graphicx}
\usepackage{float}

\newcommand{\argmin}{\operatornamewithlimits{arg\,min}}
\newcommand{\argmax}{\operatornamewithlimits{arg\,max}}

\newcommand{\demandhum}{\lambda^\text{h}}
\newcommand{\demandaut}{\lambda^\text{a}}

\newcommand{\pathidx}{i}
\newcommand{\numpaths}{N}
\newcommand{\pathset}{{[}\numpaths{]}}

\newcommand{\useridx}{j}

\newcommand{\congprof}{s}	
\newcommand{\congprofvec}{\boldsymbol{s}}	

\newcommand{\ffvelocity}{\bar{v}}
\newcommand{\numlanes}{b}
\newcommand{\spacehuman}{h^\text{h}}
\newcommand{\spaceaut}{h^\text{a}}

\newcommand{\denshum}{n^\text{h}}
\newcommand{\densaut}{n^\text{a}}
\newcommand{\flow}{f}
\newcommand{\flowhum}{f^\text{h}}
\newcommand{\flowaut}{f^\text{a}}
\newcommand{\flowhumvec}{\boldsymbol{f}^\text{h}}
\newcommand{\flowautvec}{\boldsymbol{f}^\text{a}}

\newcommand{\jamden}{\bar{n}}
\newcommand{\critdens}{\tilde{n}}
\newcommand{\capacity}{\bar{F}}

\newcommand{\autlev}{\alpha}

\newcommand{\latency}{\ell}
\newcommand{\latencyvec}{\boldell}

\newcommand{\fflatency}{a}

\newcommand{\price}{p}
\newcommand{\pricevec}{\boldsymbol{p}}
\newcommand{\dominatedroads}{D}
\newcommand{\nondominatedroads}{\bar{\dominatedroads}}
\newcommand{\reward}{r}

\newcommand{\dataset}{\mathcal{D}}
\newcommand{\flowautnocontrol}{f^\text{b}}
\newcommand{\flowautnocontrolvec}{\boldsymbol{f}^\text{b}}

\newcommand{\bp}{\boldsymbol{p}}
\newcommand{\boldell}{\boldsymbol{\ell}}

\DeclarePairedDelimiter\norm{\lVert}{\rVert}%
\newcommand\compactdots{\hbox to 1em{.\hss.\hss.}}
\newcommand{\asymeq}{\stackrel{\boldsymbol{\cdot}}{=}}

\usepackage{xcolor}
\newcommand{\ignore}[1]{}

\title{\LARGE \bf
	The Green Choice: \\ Learning and Influencing Human Decisions on Shared Roads
}

\author{Erdem B\i y\i k$^{1}$, Daniel A. Lazar$^{2}$, Dorsa Sadigh$^{1,3}$, and Ramtin Pedarsani$^{2}$
	\thanks{$^{1}$Department of Electrical Engineering, 
		Stanford University}%
	\thanks{$^{2}$Department of Electrical and Computer Engineering, 
		UC Santa Barbara}%
	\thanks{$^{3}$Department of Computer Science, 
		Stanford University}%
	\thanks{{\tt\small ebiyik@stanford.edu, dlazar@ece.ucsb.edu, dorsa@cs.stanford.edu, ramtin@ece.ucsb.edu}}
}

\begin{document}

\maketitle
\thispagestyle{empty}
\pagestyle{empty}

\begin{abstract}

Autonomous vehicles have the potential to increase the capacity of roads via platooning, even when human drivers and autonomous vehicles share roads. However, when users of a road network choose their routes selfishly, the resulting traffic configuration may be very inefficient. Because of this, we consider how to influence human decisions so as to decrease congestion on these roads. We consider a network of parallel roads with two modes of transportation: (i) human drivers who will choose the quickest route available to them, and (ii) ride hailing service which provides an array of autonomous vehicle ride options, each with different prices, to users. In this work, we seek to design these prices so that when autonomous service users choose from these options and human drivers selfishly choose their resulting routes, road usage is maximized and transit delay is minimized. To do so, we formalize a model of how autonomous service users make choices between routes with different price/delay values. Developing a preference-based algorithm to learn the preferences of the users, and using a vehicle flow model related to the Fundamental Diagram of Traffic, we formulate a planning optimization to maximize a social objective and demonstrate the benefit of the proposed routing and learning scheme.

\end{abstract}

\section{Introduction}
Road congestion is a major and growing source of inefficiency, costing drivers in the United States billions of dollars in wasted time and fuel \cite{schrank2015}. The emergence of autonomous vehicles promises to improve the efficiency of road usage to some extent by forming platoons and smoothing traffic flow. However, if drivers make their routing decisions \emph{selfishly} and seek to minimize their experienced delay, this can result in suboptimal network performance \cite{koutsoupias:1999fs}. Moreover, in the presence of selfish users, the increase in road capacity induced by replacing some human-driven vehicles with autonomous ones can paradoxically worsen average transit user delay \cite{mehr2018can}. The adverse effects arise from the way humans make decisions about how to route themselves, and the resulting inefficiency can be very large, even unbounded \cite{biyik2018altruistic}.

It is therefore important to consider the effect that self-interested users will have on the system. Previous works have shown that if autonomous users are \emph{altruistic}, average latency can be greatly reduced \cite{biyik2018altruistic}. However, it's not realistic to assume altruism in many scenarios; in this paper we consider how to monetarily incentivize such prosocial behavior. We consider the setting in which use of autonomous vehicles is offered by some benevolent ride hailing service or social planner which chooses prices for each of the routes. The users of the autonomous service will choose their routes, as well as whether or not they even want to travel, based on their time-money tradeoff. Then the remaining population of human drivers will choose routes that minimize their delay. The role of the social planner is then to choose prices that will optimize some social objective, which includes maximizing road usage and minimizing average delay. Our model can be viewed as an indirect Stackelberg game -- a game in which a social planner controls some fraction of the population's actions and the remainder of the population responds selfishly \cite{swamy2012effectiveness, krichene2017stackelberg}. However, in our model the planner can only control its portion of the vehicle flow via pricing.

In order to do so effectively, the social planner needs a model for how people make decisions between routes with various prices and latencies. We model these choices as being made based on a modification of Luce's choice axiom \cite{luce2012individual,ben1985discrete} and we model human drivers as selfish agents who reach a \emph{Wardrop Equilibrium}, a configuration in which no one could decrease their travel time by switching paths \cite{correa2011wardrop}. 

Moreover, we use a method of actively learning the choices of autonomous service users via a series of queries \cite{sadigh2017active,biyik2018batch}. This enables the planner to predict how autonomous service users will react to a set of options after a relatively low number of training queries. We verify that this method accurately predicts human choice, and that our planning algorithm can use this to improve network performance, via experiment.

We summarize our contributions as follows.
\begin{itemize}
	\item \emph{Learning Human Decisions:} We develop an active preference-based algorithm to learn how people value the tradeoff between time and money, and choose their preferred transportation option. This enables learning a model of humans' routing choices in a data-efficient manner. 
	\item \emph{Influencing Human Decisions:} We formulate and solve an optimization for ride hailing service to minimize congestion and maximize the road network utilization based on the learned model of human choice, while constraining a minimum profit for the service supplier.
	\item We validate the learning algorithm and the optimization formulation via simulations and user studies. Our results suggest carefully designed pricing schemes can provide significant improvements for traffic throughput.
\end{itemize}

\noindent \textbf{Related Work. }Previous works have shown the potential benefits that autonomous vehicles can have for traffic networks, by increasing road capacity through platooning \cite{lioris2017platoons, askari2017effect}, damping shockwaves of slowing vehicles \cite{stern2018dissipation, bhadani2018dissipation, wu2018stabilizing}, managing merges \cite{jin2018modeling}, and decongesting freeways in the event of accidents \cite{sivaranjani2015localization}. In \cite{mehr2018can}, the authors show that the capacity improvement from platooning can actually \emph{worsen} the total delay experienced by users of the networks. Relatedly, many works (\emph{e.g.} \cite{roughgarden2002bad,correa2008geometric}, on shared roads, \cite{lazar2018routing}) analyze and bound the inefficiency that can arise from network users choosing their routes selfishly instead of following some optimal routing. 

As we use financial incentives to influence this behavior, our formulation is related to work in tolling \cite{beckmann1956studies}, including works that consider users with different sensitivities to tolls \cite{fleischer2004tolls,brown2017robustness}. In contrast to many of these works, as well as some works that develop efficient algorithms for rebalancing an autonomous fleet \cite{salazar2019congestion}, we consider a probabilistic model for human choice that has been empirically shown to reflect human choice making processes \cite{daw2006cortical}, and we model vehicle flow on shared roads based on the foundational Fundamental Diagram of Traffic (FDT) \cite{daganzo1994cell,biyik2018altruistic}. This model captures the bifurcated nature of vehicle flow -- when density is low, the road is in free-flow and flow increases with vehicle density, and when density is high, the road is congested and flow decreases as density increases.

On the human choice models, there have recently been a lot of effort on learning human reward functions which are assumed to be sufficient to model the preferences. Inverse reinforcement learning \cite{abbeel2004apprenticeship,ng2000algorithms,abbeel2005exploration,ziebart2008maximum} and preference-based learning \cite{sadigh2017active,biyik2018batch,akrour2012april,christiano2017deep} are the most popular choices. In this paper, we employ preference-based learning as it is a natural fit to our problem. We actively synthesize queries as a non-trivial generalization and extension of \cite{sadigh2017active} for data-efficiency and better usability.

\section{Problem Setting}
\begin{figure}
	\centering
	\includegraphics[width=\linewidth]{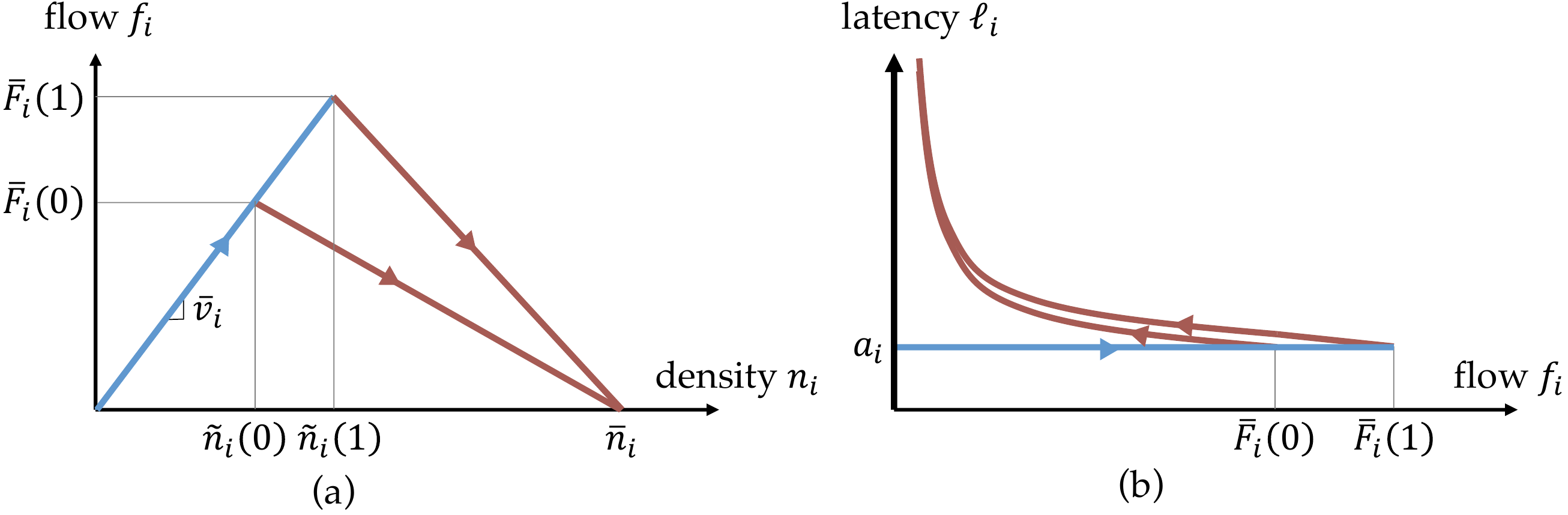}
	\vspace{-15px}
	\caption{Vehicle flow model for mixed autonomy based on the Fundamental Diagram of Traffic. (a) and (b) show flow vs density and latency vs flow characteristics if all vehicles are either human-driven or autonomous, where autonomy increases the critical density $\critdens_\pathidx$ and the maximum flow $\capacity_\pathidx$ of a road. The blue and red regions correspond to the road being free-flow and congested states, respectively, and the arrow denotes increasing density.}
	\label{fig:FDTs}
\end{figure}

\noindent \textbf{Vehicle Flow Model. }We consider a nonatomic model in which individual users control an infinitesimal unit of vehicle flow and this quantity is treated as a continuous variable. To model flow of mixed-autonomous traffic on roads, we use our model in \cite{biyik2018altruistic}, summarized in this section. This is an extension of the Fundamental Diagram of Traffic \cite{daganzo1994cell}, in which flow inhabits one of two regimes -- free flow, in which the flow on a road \emph{increases} as vehicle density increases, and congestion, in which flow \emph{decreases} as density increases. The flow decreases until the \emph{jam density}, $\jamden_\pathidx$, at which point the flow ceases entirely. Similarly, the relationship between latency and flow on a road shows a similar divide: in free flow latency is constant with respect to flow, and in congestion latency \emph{increases} as flow decreases. These behaviors are shown in Fig.~\ref{fig:FDTs}. $\ffvelocity_\pathidx$ denotes the free-flow speed of the road.

\begin{figure}
	\centering
	\vspace*{-10px}
	\includegraphics[width=0.75\linewidth]{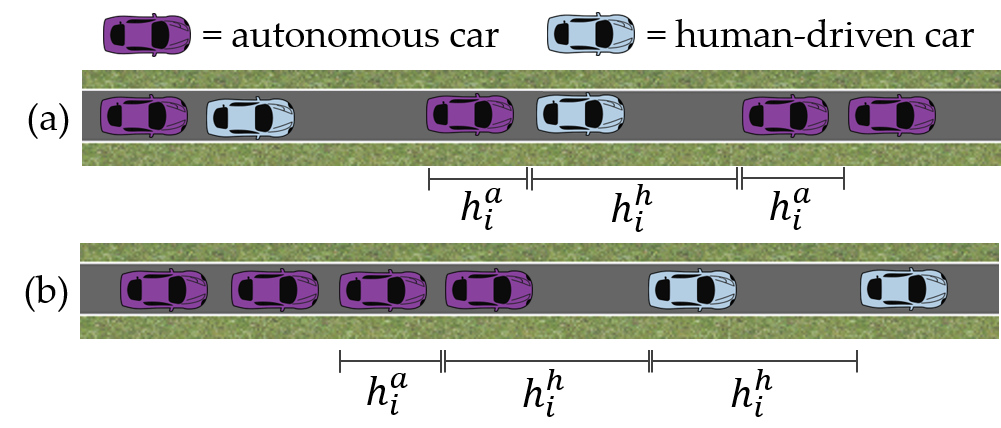}
	\vspace*{-5px}
	\caption{Illustration of vehicle spacing if (a) autonomous vehicles maintain the same headway behind any vehicle they follow, regardless of type, or (b) autonomous vehicles can reorder themselves to form a large platoon. The headway that the vehicles keep informs the critical density of a road, or the density at which a road becomes congested.}
	\vspace*{-15px}
	\label{fig:spacing}
\end{figure}

The transition between free flow and congestion occurs when the density on a road exceeds the \emph{critical density} of that road, denoted $\critdens_\pathidx$. We model the effect of autonomous vehicles as increasing this critical density, since autonomous cars can maintain shorter headways than human drivers, as follows.
\begin{equation}\label{eq:crit_dens}
\critdens_{\pathidx}(\autlev_{\pathidx}) := \frac{ \numlanes_{\pathidx}}{\autlev_{\pathidx}\spaceaut_\pathidx + (1-\autlev_{\pathidx})\spacehuman_\pathidx} \; ,
\end{equation}
where $\autlev_{\pathidx} \in {[0,1{]}}$ is the \emph{autonomy level}, or fraction of vehicles on road $\pathidx$ that are autonomous, $\numlanes_{\pathidx}$ is the number of lanes on that road, and $\spaceaut_\pathidx$ and $\spacehuman_\pathidx$ are the space occupied by an autonomous and human driven car, respectively, when traveling at the nominal velocity of the road. This model \cite{lazar2018routing} assumes that either the autonomous vehicles maintain the same headway to any vehicle they may follow or that they can rearrange themselves on a road to form large platoons \cite{lazar2018price,lazar2018maximizing}, as shown in Fig.~\ref{fig:spacing}. Note that if autonomous vehicles can keep shorter headways than human drivers (\emph{i.e.} $\spaceaut_\pathidx < \spacehuman_\pathidx$), the presence of autonomous cars increases the critical density and therefore also increases the maximum flow possible on a road, i.e. $\capacity_\pathidx(1) > \capacity_\pathidx(0)$ where $\capacity_\pathidx$ is a function of autonomy level that gives the maximum capacity of road $\pathidx$. Fig.~\ref{fig:FDTs} shows the effect of increased autonomy on a road.

We use $\congprof_\pathidx$ to denote the state of road $i$ such that $s_i=0$ denotes it is in free flow and $s_i=1$ denotes it is congested. We use $\denshum_\pathidx$ and $\densaut_\pathidx$ to denote the density of human-driven and autonomous vehicles, respectively on road $\pathidx$. Accordingly,
\begin{equation*}
\congprof_\pathidx = \begin{cases}
0 & \denshum_\pathidx + \densaut_\pathidx \le \critdens_\pathidx(\densaut_\pathidx/(\denshum_\pathidx+ \densaut_\pathidx)) \\
1 & \text{otherwise} \; .
\end{cases}
\end{equation*}
We model the road as first-in-first-out, so $\autlev_{\pathidx} = \densaut_\pathidx/(\denshum_\pathidx + \densaut_\pathidx)= \flowaut_\pathidx/(\flowhum_\pathidx + \flowaut_\pathidx)$, where $\flowhum_\pathidx$ and $\flowaut_\pathidx$ respectively are the human-driven and autonomous vehicle flow on road $\pathidx$.

The flow on a road, $\flow_\pathidx = \flowhum_\pathidx + \flowaut_\pathidx$ is a function of the density of each vehicle type as follows.
\begin{align}\label{eq:flow}
&\flow_\pathidx(\denshum_\pathidx, \densaut_\pathidx) := \nonumber \\
&\begin{cases}
\ffvelocity_\pathidx \cdot (\denshum_\pathidx + \densaut_\pathidx) ,& \text{if } \denshum_\pathidx + \densaut_\pathidx \le \critdens_{\pathidx}(\autlev_{\pathidx})\\
\frac{\ffvelocity_\pathidx \cdot \critdens_{\pathidx}(\autlev_{\pathidx}) \cdot (\jamden_\pathidx - (\denshum_\pathidx + \densaut_\pathidx)) }{\jamden_\pathidx - \critdens_{\pathidx}(\autlev_{\pathidx})} ,              & \text{if }  \critdens_{\pathidx}(\autlev_{\pathidx}) \le \denshum_\pathidx + \densaut_\pathidx \le \jamden_\pathidx \\
0 ,& \text{otherwise } \; .
\end{cases}
\end{align}

We can then write the latency as a function of vehicle flow and the state of the road \cite{krichene2017stackelberg,biyik2018altruistic}:
\begin{align*}
&\latency_\pathidx(\flowhum_\pathidx,\flowaut_\pathidx,\congprof_\pathidx) = \begin{cases}
\frac{d_\pathidx}{\ffvelocity_\pathidx} & \text{if } \congprof_\pathidx = 0 \\
d_\pathidx\left(\frac{\jamden_\pathidx}{\flowhum_\pathidx + \flowaut_\pathidx} + \frac{\critdens_\pathidx(\autlev_{\pathidx}) - \jamden_\pathidx}{\ffvelocity_\pathidx \cdot \critdens_\pathidx(\autlev_{\pathidx})}\right) & \text{if } \congprof_\pathidx = 1 \; .
\end{cases}
\end{align*}
where $d_\pathidx$ is the length of road $\pathidx$.

We now know how to relate the vehicle flow to road latency so we can start reasoning about how people choose their routes.

\noindent \textbf{Problem Overview. }We describe the various facets of the problem, preceded by a high-level overview of the objective. We assume that the demand of human drivers is fixed, and that the demand of people using the autonomous service is elastic -- \emph{if prices and latencies are high, some people may choose not to use the autonomous mobility service.} Our goal is then to simultaneously maximize the number of autonomous service users that can use the road, and minimize the average delay experienced by all the people using the roads.

Our control variables in this optimization are the \emph{latencies} on each road and the \emph{prices} offered to the users for traveling on the routes. However, we can't arbitrarily choose prices and latencies -- we need to respect 
\begin{enumerate}
	\item \emph{the characteristics of the roads}, in terms of how the flow demand for a road corresponds to the latency on that road, and 
	\item \emph{how people make decisions}, making sure that the number of people who choose each option corresponds to the latencies of the roads described in the options.
\end{enumerate}
 Moreover, we want to be fair and must therefore offer the same pricing and routing options to all users of the autonomous service. 

How a population of users chooses between a variety of price and latency pairs depends on their valuation of time and money. Without knowing this choice model we cannot plan vehicle flows and will not be able to ensure the resulting configuration matches our vehicle flow models for the roads. Also, since different populations may have different parameter distributions, we need to learn this tradeoff for our population before we can estimate how many people will choose which option. To untangle these constraints, in the subsequent sections we describe our models for vehicle flow and human choices. Following that we return with a mathematical formulation of this optimization.

\section{Learning the Human Choice Model}
\label{sec:choice_model}

In this section we describe how both human drivers and the autonomous service users make routing decisions. We then provide a data-efficient learning algorithm for estimating the way that the autonomous service users make routing decisions, allowing us to predict the response a population will have to a set of latency/price options.

\begin{figure}
	\centering
	\includegraphics[width=0.7\linewidth]{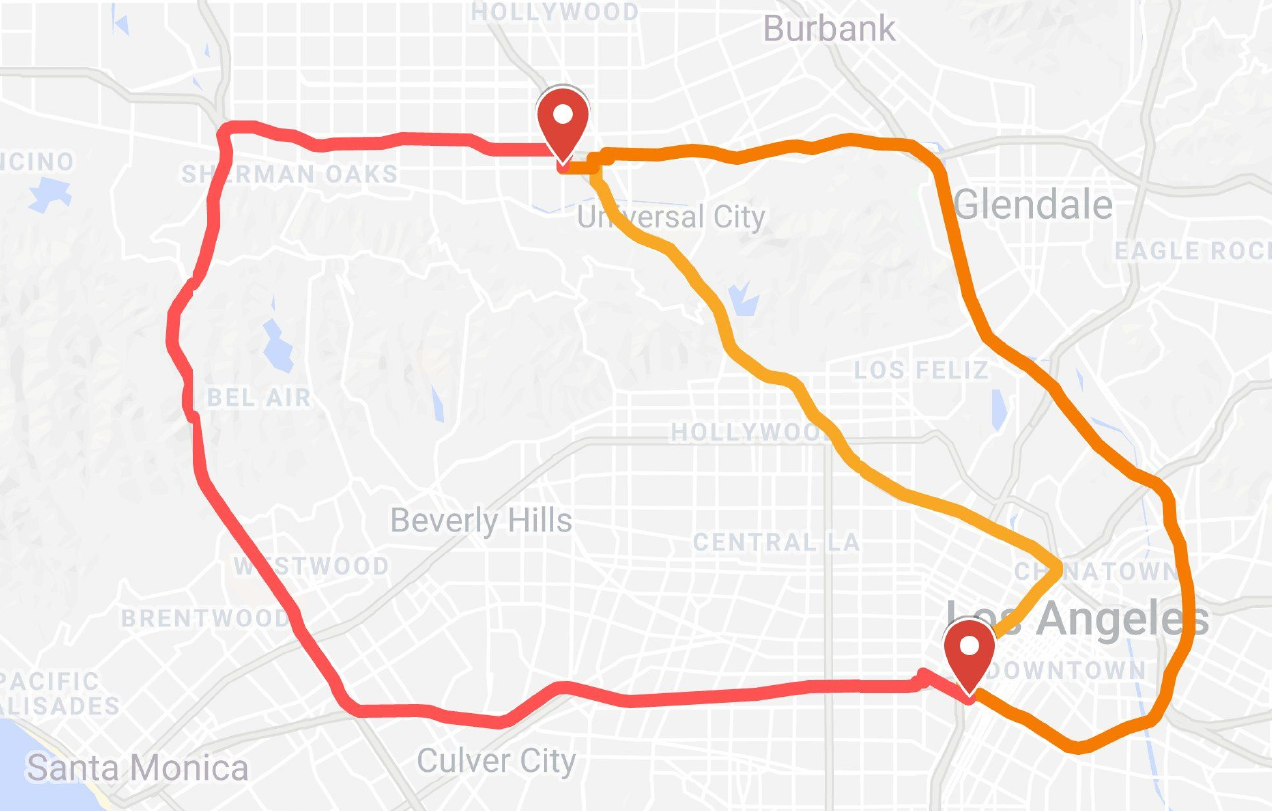}
	\caption{An example, from Los Angeles, of a network of parallel roads.}
	\label{fig:network}
\end{figure}

\noindent \emph{Human Drivers. }We consider a network of $\numpaths$ parallel roads, an example of which is shown in Fig.~\ref{fig:network}.
We assume that no two roads have the same free-flow latency and denote the free-flow latency of path $\pathidx$ by $\fflatency_\pathidx$. We order the indices such that $\fflatency_1 < \fflatency_2 < \ldots < \fflatency_\numpaths$, and we use $[k]$ to denote the set of the first $k$ roads.

\begin{figure}
	\centering
	\includegraphics[width=\linewidth]{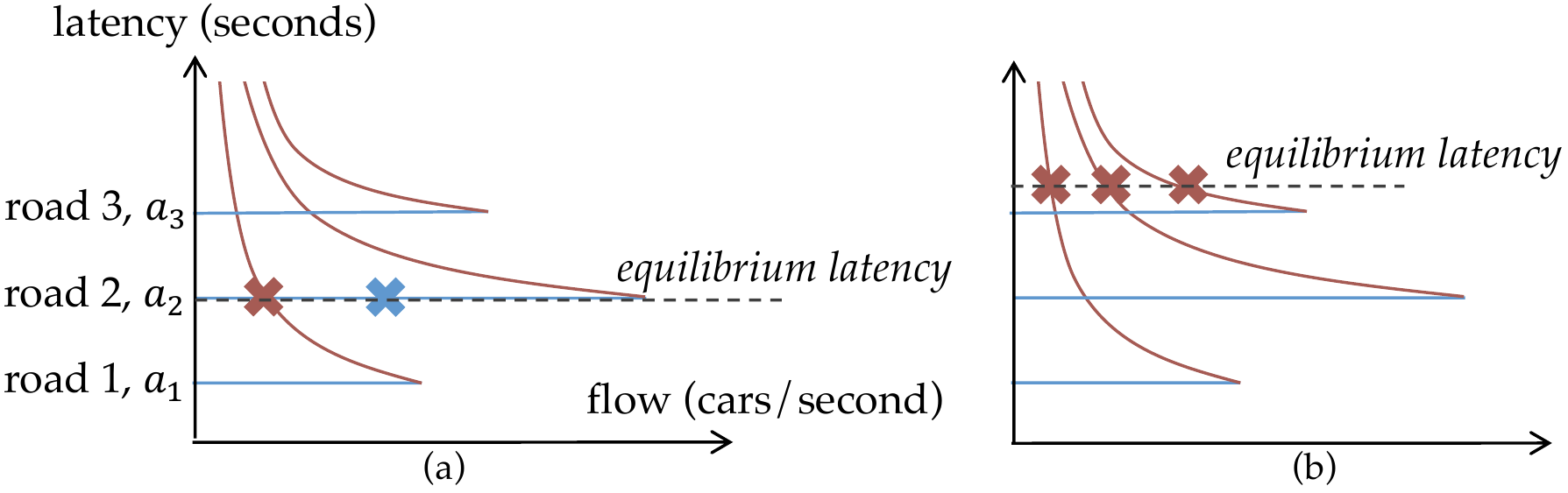}
	\caption{Some possible equilibria of a three-road network with fixed flow demand. Blue and red lines denote the free-flow and congested regimes, respectively. Equilibria may have (a) one road in free-flow or (b) all used roads may be congested. An equilibrium has an associated \emph{equilibrium latency} experienced by all selfish users. By considering a given equilibrium delay, we can reason about which roads must be congested at that equilibrium as well how much flow is on each road.}
	\vspace*{-15px}
	\label{fig:equilibria}
\end{figure}

For the human drivers, we assume their only consideration is minimizing their commute time. This leads to a \emph{Wardrop Equilibrium} \cite{wardrop1900some,depalma1998optimization}, which, on parallel roads, means that if one road has positive flow on it, all other roads must have higher or equal latency. Formally, 
$$\flowhum_\pathidx>0 \implies \latency_\pathidx(\flowhum_\pathidx, \flowaut_\pathidx, \congprof_\pathidx) \le \latency_{\pathidx'}(\flowhum_{\pathidx'}, \flowaut_{\pathidx'}, \congprof_{\pathidx'}) \; \forall \pathidx,\pathidx' \in \pathset \; .$$
This implies that all selfish users experience the same latency. It is therefore useful to consider the flow-latency diagrams of roads when studying which equilibria exist -- by fixing the latency on the y-axis, one can reason about which roads must be congested to achieve that equilibrium. As shown in Fig.~\ref{fig:equilibria}, equilibria may have one road in free-flow and rest congested, or all may be congested, assuming no two roads have the same free-flow latency \cite{krichene2017stackelberg, biyik2018altruistic}.

\noindent \emph{Autonomous Service Users. }
Though human drivers are motivated directly only by delay, users of the autonomous service experience cost in both delay and the price of the ride. We model the users as having some underlying reward function, which is parameterized by their time/money tradeoff, as well as the desirability of the option of traveling by some other means such as walking or public transit. We assume that a strictly dominated option in terms of latency and prices is completely undesirable. Formally, this set of dominated roads is defined as 
\begin{align*}
D = \{\pathidx \in \pathset \: |\:(\price_\pathidx>\price_{\pathidx'} \land \latency_\pathidx \geq \latency_{\pathidx'}) \lor (\price_\pathidx \geq \price_{\pathidx'} \land \latency_\pathidx > \latency_{\pathidx'})\\ \textrm{ for some } \pathidx'\in \pathset\} \; .
\end{align*}
We also define the complementary set, i.e., nondominated roads, $\nondominatedroads = \pathset \backslash \dominatedroads$. We model the reward function of user $\useridx$ for choosing road $\pathidx$ as follows:
\begin{align*}
&\reward_\useridx(\boldsymbol{\latency},\mathbf{\price},\pathidx) = \\
\quad &\begin{cases}
- {\omega_\useridx}_1 \latency_\pathidx - {\omega_\useridx}_2 \price_\pathidx & \textrm{if user }\useridx \textrm{ chooses road }\pathidx\textrm{ and }\pathidx \in \nondominatedroads \\
-\infty & \textrm{if user }j\textrm{ chooses road }\pathidx\textrm{ and }\pathidx \in \dominatedroads\\
-\zeta_\useridx \latency^w & \textrm{if user }\useridx \textrm{ declines the service ,}
\end{cases}
\end{align*}
where $\latencyvec$ denotes the vector of road latencies, $\pricevec$ denotes prices, ${\boldsymbol{\omega}_\useridx} = \begin{bmatrix}{\omega_\useridx}_1 & {\omega_\useridx}_2 \end{bmatrix}^T$ characterizes the users' time/money tradeoff and $\zeta_\useridx$ specifies their willingness to use an alternative option, which has delay $\latency^w$. This alternate option could be walking, biking, or public transportation; the nature of this option determines the delay as well and influences $\zeta_\useridx$.

We don't assume that users are simple reward maximizers. Rather, we draw from Luce's Choice Axiom \cite{luce2012individual} to model the probability with which users choose each option, in combination with a softmax rule \cite{daw2006cortical}, formalized as follows.
\begin{equation}\label{eq:decision_probs}
P(\textrm{user }\useridx \textrm{ chooses route } \pathidx) = \frac{\exp\left(\reward_\useridx(\latencyvec,\pricevec,\pathidx)\right)}{\sum_{\pathidx' \in \pathset} \exp\left(\reward_\useridx(\latencyvec,\pricevec,\pathidx')\right)} \; .
\end{equation}

In order to determine optimal pricing, we want to know how many users will choose each route as a function of the route prices and latencies. To this end, we define $q_i(\latencyvec, \pricevec)$ as the expected fraction of autonomous service users that will choose route $\pathidx$. If the parameter distribution for autonomous service users is $g(\boldsymbol{\omega},\zeta)$, then
\begin{align*}
&q_i(\latencyvec, \pricevec) = \int_0^\infty\int_0^\infty\int_0^\infty \\
& \; g(\boldsymbol{\omega},\zeta)P(\textrm{user chooses route } \pathidx | \boldsymbol{\omega}, \zeta) d{\omega_1} d{\omega_2} d\zeta \; .
\end{align*}

This expression relate prices and latency to human choices, enabling us to determine the prices that will maximize our social objective. This will be important in constraining our optimization to only consider latency/price options that correspond to the desired vehicle flows.

\noindent \textbf{Data-Efficient Learning of Human Reward Functions. }
While the routes in a specific network can be fully modeled with the physical properties and the speed limits, user's decision models must be learned in a data-driven way. As described above, reward function parameters $(\boldsymbol{\omega}, \zeta)$ are sufficient to model the decisions of a user. 

The parameters will be different among the users. While a business executive might prefer paying extra to reach their company quickly, a graduate student may decide to go to the lab a little later in order to save a few dollars. Therefore, we have to learn personalized parameters $\boldsymbol{\omega}$ and $\zeta$. Although it might be possible to make a rough guess based on personal information, such as age, income, occupation, location, etc; there are a few problems about this approach: 1) people might be unwilling to share those information accurately, 2) the variance will be high, since the decisions are also affected by the personal traits that are hard to quantify, 3) we do not know how we will instantiate the parameters, because putting this as a regression problem would require some users to explicitly quantify their $\boldsymbol{\omega}$ and $\zeta$, which is rather unrealistic.

Due to these problems, we choose to learn the parameters from users' previous choices, which is known as \emph{preference-based learning}. If user $\useridx$ chooses from a variety of options, the user's choice gives us a noisy estimate of which road $\pathidx \in \pathset$ maximizes their reward function $\reward_\useridx(\boldsymbol{\latency}, \mathbf{\price}, \pathidx)$. We could start from either uniform priors or priors informed by domain knowledge, then sample from the distribution $g(\boldsymbol{\omega},\zeta)$.

However, a major drawback of doing so is how quickly we learn the user preferences. Preference-based learning suffers from the small amount of information that each query carries. For example, if we show $4$ options to a user (including the option to decline the service), then the maximum information we can get from that query is just $2$ bits. In fact, the information gained can be much lower in later queries, because the posterior is already a better estimate than the original prior distribution. To tackle this issue, previous works pose the query generation/synthesis problem as a submodular optimization and maximize a measure of the expected change in the learned distribution after the queries. \cite{sadigh2017active,biyik2018batch}. 

While those works focused only on pairwise queries, in this case we expect to pose several route options to the users and therefore need more general query types. By using these more general queries that offer a variety of routes with varied latencies and prices, we can consider a number of ways of using this learning framework to learn the human preference distribution.
\begin{itemize}
	\item We could do a user study on a few people to learn a good prior distribution. While we could achieve this using random queries, it would take too many queries to learn a useful distribution and most of the queries would not be sensible. By utilizing active learning, we are able to learn the parameters of each user using only a few queries.
	\item We could do an exploration/exploitation strategy if we are allowed to break the fairness constraint a few times for some small portion of the users. This could be possible through user-specific campaigns. That is, for each user we may either choose to use the learned model or to offer special rates that would help us profile the user better.
	\item We could do an initial profiling study for each new user, where we ask them queries synthesized by the active learning framework.
\end{itemize}

To implement any of these options, we formulate the following active learning optimization. First, we discuss the general preference-based learning framework. Given the data from previous choices of user $\useridx$, which we denote as $\dataset_\useridx$, we can formalize the probability of $(\boldsymbol{\omega}_j,\zeta_j)$ being true parameters for that user as follows:
\begin{align*}
P(\boldsymbol{\omega}_\useridx,\zeta_\useridx | \dataset_\useridx) &\propto P(\boldsymbol{\omega}_\useridx,\zeta_\useridx)P(\dataset_\useridx | \boldsymbol{\omega}_\useridx,\zeta_\useridx)\\
&= P(\boldsymbol{\omega}_\useridx,\zeta_\useridx)\prod_{m}P({\dataset_\useridx}_m | \boldsymbol{\omega}_\useridx,\zeta_\useridx)
\end{align*}
where ${\dataset_\useridx}_m$ denotes the route user $\useridx$ chose in their $m^{\textrm{th}}$ choice (with ${\dataset_\useridx}_m=0$ meaning that the user saw the options, but preferred the alternative option). The first expression is due to Bayes' rule and the equality is due to the assumption that the users' choices are conditionally independent from each other given the reward function parameters.

We can readily obtain the second term from the human choice model. For the prior, we can simply use a uniform distribution over nonnegative values of the parameters. The prior can be crucial especially when we do not have enough data for a new user. In such settings, we incorporate domain knowledge to start with a better prior than a uniform distribution.

We can then use this unnormalized $P(\boldsymbol{\omega}_\useridx,\zeta_\useridx | \dataset_\useridx)$ to obtain the samples of $(\boldsymbol{\omega}_\useridx,\zeta_\useridx)$ using Metropolis-Hastings algorithm. Doing this for each user, which can be easily parallelized, we can directly obtain the samples $(\boldsymbol{\bar\omega},\bar\zeta)\sim g(\boldsymbol{\omega},\zeta)$.

Next we formulate the active learning framework, which is needed so that it will not take an excessive number of queries to learn human preferences. For this, we want to maximize the expectation of the difference between the prior and the unnormalized posterior:
\begin{align*}
\textrm{query}_m^* &= \argmax_{\textrm{query}_m} \mathbb{E}_{{\dataset_\useridx}_m}\Big[P(\boldsymbol{\omega}_\useridx,\zeta_\useridx|{\dataset_\useridx}_{1:m-1})-\\
&\qquad P(\boldsymbol{\omega}_\useridx,\zeta_\useridx|{\dataset_\useridx}_{1:m-1})P({\dataset_\useridx}_m | \boldsymbol{\omega}_\useridx,\zeta_\useridx,{\dataset_\useridx}_{1:m-1})\Big]\\
&= \argmin_{\textrm{query}_m} \mathbb{E}_{{\dataset_\useridx}_m}\left[P({\dataset_\useridx}_m | \boldsymbol{\omega}_\useridx,\zeta_\useridx,{\dataset_\useridx}_{1:m-1})\right]
\end{align*}
As we will use the sampled $(\boldsymbol{\omega}_\useridx,\zeta_\useridx)$ for computing the probabilities of each route choice for the expectation over ${\dataset_\useridx}_m$, we can write the optimization as:
\begin{align*}
\textrm{query}_m^* &\asymeq \argmin_{\textrm{query}_m} \mathbb{E}_{{\dataset_\useridx}_m}\!\left[\sum_{\boldsymbol{\bar\omega}_\useridx,\bar\zeta_\useridx}\!P({\dataset_\useridx}_m | \boldsymbol{\bar\omega}_\useridx,\bar\zeta_\useridx,{\dataset_\useridx}_{1:m\!-\!1})\right]
\end{align*}
where we have $M$ samples denoted as $(\boldsymbol{\bar\omega}_\useridx,\bar\zeta_\useridx)$, $\asymeq$ denotes asymptotic equality as $M\to\infty$, and the term $1/M$ is canceled. Using the law of total probability, we can also write
\begin{align*}
P({\dataset\!_\useridx}_m& | {\dataset\!_\useridx}_{1:m-1}) \\&= \sum_{\boldsymbol{\omega}_\useridx,\zeta_\useridx}P(\boldsymbol{\omega}_\useridx,\zeta_\useridx|{\dataset\!_\useridx}_{1:m-1})P({\dataset\!_\useridx}_m | {\dataset_\useridx}_{1:m-1},\boldsymbol{\omega}_\useridx,\zeta_\useridx)\\
&\asymeq\frac1M\sum_{\boldsymbol{\bar\omega}_\useridx,\bar\zeta_\useridx}P({\dataset_\useridx}_m | \boldsymbol{\bar\omega}_\useridx,\bar\zeta_\useridx,{\dataset_\useridx}_{1:m-1})
\end{align*}
which then leads to
\begin{align*}
&\textrm{query}_m^*\\
&\!\asymeq\!\argmin_{\textrm{query}_m}\! \sum_{{\dataset\!_\useridx}_m}P({\dataset\!_\useridx}_m|{\dataset\!_\useridx}_{1:m\!-\!1})\!\sum_{\boldsymbol{\bar\omega}_\useridx,\bar\zeta_\useridx}\!P({\dataset\!_\useridx}_m | \boldsymbol{\bar\omega}_\useridx,\bar\zeta_\useridx,{\dataset\!_\useridx}_{1:m\!-\!1})\\
&\!\asymeq\!\argmin_{\textrm{query}_m}\! \sum_{{\dataset_\useridx}_m}\left(\sum_{\boldsymbol{\bar\omega}_\useridx,\bar\zeta_\useridx}P({\dataset_\useridx}_m | \boldsymbol{\bar\omega}_\useridx,\bar\zeta_\useridx)\right)^2
\end{align*}
We can easily compute this objective value for any given $\textrm{query}_m$. This optimization is unfortunately nonconvex due to the human choice model. As in previous works, we assume local optima is good enough \cite{sadigh2017active, biyik2018batch}. We then use a Quasi-Newton method (L-BFGS \cite{andrew2007scalable}) to find the local optima, and we repeat this for $1000$ times starting from random initial points so that we can get closer to the global optimum.

\section{Influencing Human Decision Making}

\noindent \textbf{Problem Formulation. } We are now ready to formulate our objective. Since the demand of autonomous service users is elastic, we cannot just minimize the average latency. That would result in extremely high prices to reduce latency by keeping autonomous service users off the network. Hence, we consider the objective to be a combination of maximizing road usage and minimizing average travel time. We parameterize this tradeoff with parameter $\theta \ge 0$ in the cost function
\begin{equation*}
J(\flowhumvec,\!\flowautvec,\!\congprofvec) \!=\! \frac{\sum_{\pathidx \in \pathset}{(\flowhum_\pathidx\!+\!\flowaut_\pathidx) \latency_\pathidx(\flowhum_\pathidx, \flowaut_\pathidx, \congprof_\pathidx)}}{\sum_{\pathidx \in \pathset}{(\flowhum_\pathidx+\flowaut_\pathidx)}} - \theta\!\sum_{\pathidx \in \pathset}\! (\flowhum_\pathidx\!+\!\flowaut_\pathidx) \; .
\end{equation*}

Given $\mathbf{q}(\latencyvec,\pricevec)$, the cost function $J(\flowhumvec, \flowautvec, \congprofvec)$, inelastic flow demand of human drivers $\demandhum$, and elastic demand of autonomous users $\demandaut$, we are ready to formulate the planning optimization. The most straightforward way is to optimize jointly over $\flowhumvec,\flowautvec,\latencyvec$ and $\bp$. However, $\boldell$ is fully defined by $\flowhumvec,\flowautvec$ and $\congprofvec$. Hence, instead of $\latencyvec$, we use $\congprofvec$, which will help as the optimization is nonconvex and we can search over the possible values of $\congprofvec$. The problem is then formulated as:
\begin{flalign}
\min_{\flowhumvec,\flowautvec,\pricevec \in \mathbb{R}^\numpaths_{\ge 0},k \in \pathset,\congprof_{k:N}\in\{0,1\}^{N-k+1}}  J(\flowhumvec,\flowautvec,\congprofvec)&&
\label{eq:optimization}
\end{flalign}
\begin{align}
\textrm{subject to } & \sum_{\pathidx \in \pathset} \flowhum_\pathidx = \demandhum \label{c:1}\\
& \flowaut_\pathidx = \demandaut q_i(\latencyvec(\flowhumvec,\flowautvec,\congprofvec), \pricevec), \forall \pathidx \in \pathset \label{c:2}\\
& a_k \leq \latency_k(\flowhum_k,\flowaut_k,\congprof_k) \leq a_{k+1} \label{c:3}\\
& \flowhum_\pathidx=0, \forall \pathidx \in \pathset \setminus {[}k{]} \label{c:4}\\
& \latency_\pathidx(\flowhum_\pathidx,\flowaut_\pathidx,1) = \latency_k(\flowhum_k,\flowaut_k,\congprof_k), \forall \pathidx \in {[}k-1{]} \label{c:5}\\
& \flowhum_\pathidx + \flowaut_\pathidx \leq \capacity\left(\flowaut_\pathidx/(\flowaut_\pathidx+\flowhum_\pathidx)\right), \forall \pathidx \in \pathset \label{c:6}\\
& \sum_{\pathidx \in \pathset} \left(\flowaut_\pathidx p_\pathidx - \flowaut_\pathidx d_\pathidx c\right) \geq \bar{P} \label{c:7}
\end{align}
with $\congprof_{1:k-1}=1$ due to selfishness of human-driven vehicles, where $k$ is the longest road with human-driven vehicle flow, $\bar{P}$ is the minimum profit per unit time that we want to make from autonomous service users and $c$ is the constant fuel cost per unit length, which can be easily calculated by multiplying the vehicles' fuel consumption per unit distance and the fuel price. We can describe the constraints as follows.
\begin{enumerate}[nosep]
	\setcounter{enumi}{4}
	\item The human-driven vehicle flow demand is fulfilled.
	\item Autonomous flow will be distributed based on the choice model described in the preceding section.
	\item The ``longest equilibrium road" has latency on the given interval of free flow latencies.
	\item Human-driven cars are selfish, i.e. no human-driven car will experience higher latencies than the road $k$.
	\item The congested roads have the same latency as the ``longest equilibrium road".
	\item The maximum capacities of the roads are respected.
	\item The minimum profit per unit time is satisfied.
\end{enumerate}

We can further improve the search space by relying on the heuristic that the roads that are not used by the human-driven vehicles will be in free-flow, i.e. $s_{k+1:N}=0$. While we do not have a proof for this, we also note constructing counterexamples seems to require extremely careful design, which suggests the heuristic holds in general. Furthermore, the following proposition shows we could also set $s_k=0$ under an additional assumption.
\begin{proposition}
	Under the assumptions that the reward function is homogeneous among the users and $\omega_2>0$, there exists a free-flow road $k$ in the optimal solution to the problem such that $\latency_{\pathidx}=\latency_k$ for $\forall \pathidx\in[k]$, and $\flowhum_{\pathidx} = 0$ for $\forall \pathidx \in\pathset\setminus[k]$ as long as the optimization is feasible.
	\label{prop:eq_road}
\end{proposition}
The proof is given in Appendix.

\noindent\textbf{Generalizations.} We assumed all autonomous cars are controlled by a centralized social planner. To extend our framework to scenarios where this is not the case and the social planner has the control over a fraction of autonomous cars, we can simply do the following modifications: The optimization will also be over $\flowautnocontrolvec\in\mathbb{R}_{\geq0}^\numpaths$, which will now represent the autonomous flow that does \emph{not} belong to the planner. We add the corresponding constraint of \ref{c:1} for $\flowautnocontrol_\pathidx$. Similar to \ref{c:4}, we will also have $\flowautnocontrol_\pathidx=0$ for $\pathidx \in \pathset \setminus {[}k{]}$ due to selfishness. We will also replace $\flowaut$'s in \ref{c:3}, \ref{c:5} and \ref{c:6} with $\flowaut + \flowautnocontrol$ with appropriate subscripts. These simple modifications enable a more general use.

\noindent \textbf{Solving the Optimization. }
After learning the distribution $g(\boldsymbol{\omega},\zeta)$, we first take $M$ samples $(\boldsymbol{\bar\omega},\bar\zeta)\sim g(\boldsymbol{\omega},\zeta)$ as we previously described in Section~\ref{sec:choice_model}. Using these samples, we approximate the expected fraction of autonomous users that will choose route $\pathidx$ as:
\begin{align*}
&q_i(\latencyvec, \pricevec) \asymeq \frac{1}{M}\sum_{\boldsymbol{\bar\omega},\bar\zeta} P(\textrm{user chooses route } \pathidx | \boldsymbol{\bar\omega}, \bar\zeta)
\end{align*}

We then locally solve the nonconvex planning optimization using interior point algorithms \cite{hribar1997interior,byrd2000trust,waltz2006interior}. We execute the algorithms with $100$ random initial points for each run to get closer to global optimum.

\section{Experiments \& Results}
\label{sec:results}
To validate our framework, we conducted different simulations and a user study.

\noindent\textbf{Hypotheses.} We test three hypotheses that together suggest our framework successfully reduces traffic congestion through pricing, after it effectively learns the humans' choice models:\\
\indent\textbf{H1}: Our active learning algorithm is able to learn the presented human choice model of autonomous service users in a data-efficient way.\\
\indent\textbf{H2}: Our planning optimization reduces the latencies by creating altruistic behavior through pricing.\\
\indent\textbf{H3}: When used by humans, the overall framework works well and is advantageous over non-altruistic algorithms.

\noindent\textbf{Implementation Details.} In the planning optimization, we used the heuristic that $s_{k+1:N}=0$, but did not set $s_k=0$. We assumed there is only one alternative option for autonomous service users, and it is walking to the destination. We set $c=6\times 10^{-5}$ USD/meter, $\mathbf{b}=1$. We assumed the human-driven cars keep a $2$-second headway distance with the leading car, whereas autonomous cars keep $1$-second distance. We set the length of the cars as $5$ meters, and the minimum gap between two cars as $2$ meters.

\begin{figure}[h]
	\centering
	\vspace{-10px}
	\includegraphics[width=\linewidth]{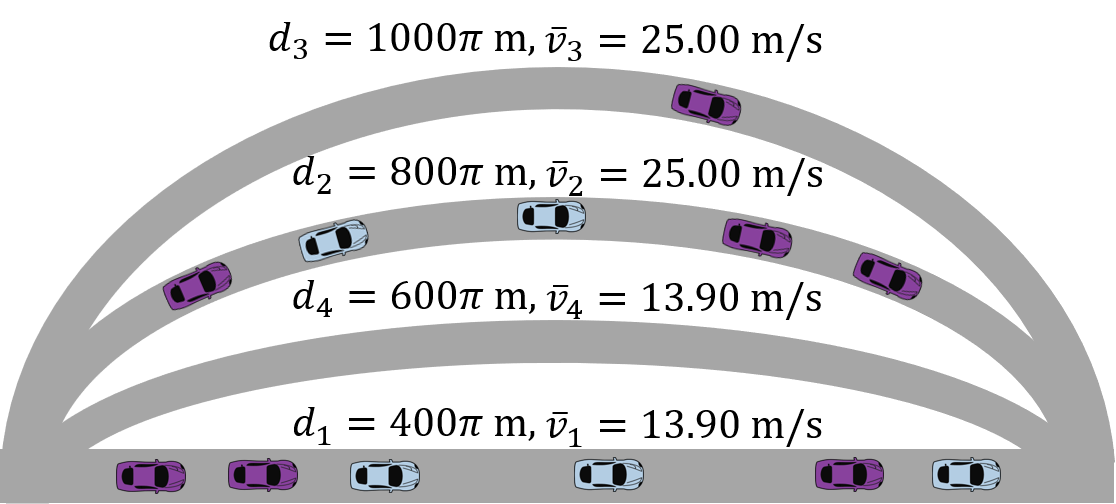}
	\caption{The $4$-road network from \cite{biyik2018altruistic}. The roads are not to the scale and ordered with respect to the free-flow latencies.}
	\vspace{-15px}
	\label{fig:wafr_network}
\end{figure}

\noindent\textbf{Experiments and Analyses.} To test \textbf{H1}, we simulated $5$ autonomous service users with different preferences. We tested our data-efficient learning framework by asking two sets of $200$ queries, each of which consisted of $4$ road options, similar to Fig.~\ref{fig:wafr_network}, and a walking option, to the simulated users. The queries were generated actively in the first set and randomly in the second. After each query, we recorded the sample $(\bar{\boldsymbol{\omega}},\bar\zeta)$ which has the highest likelihood as our estimates.

\begin{figure}[h]
	\centering
	\vspace{-10px}
	\includegraphics[width=\linewidth]{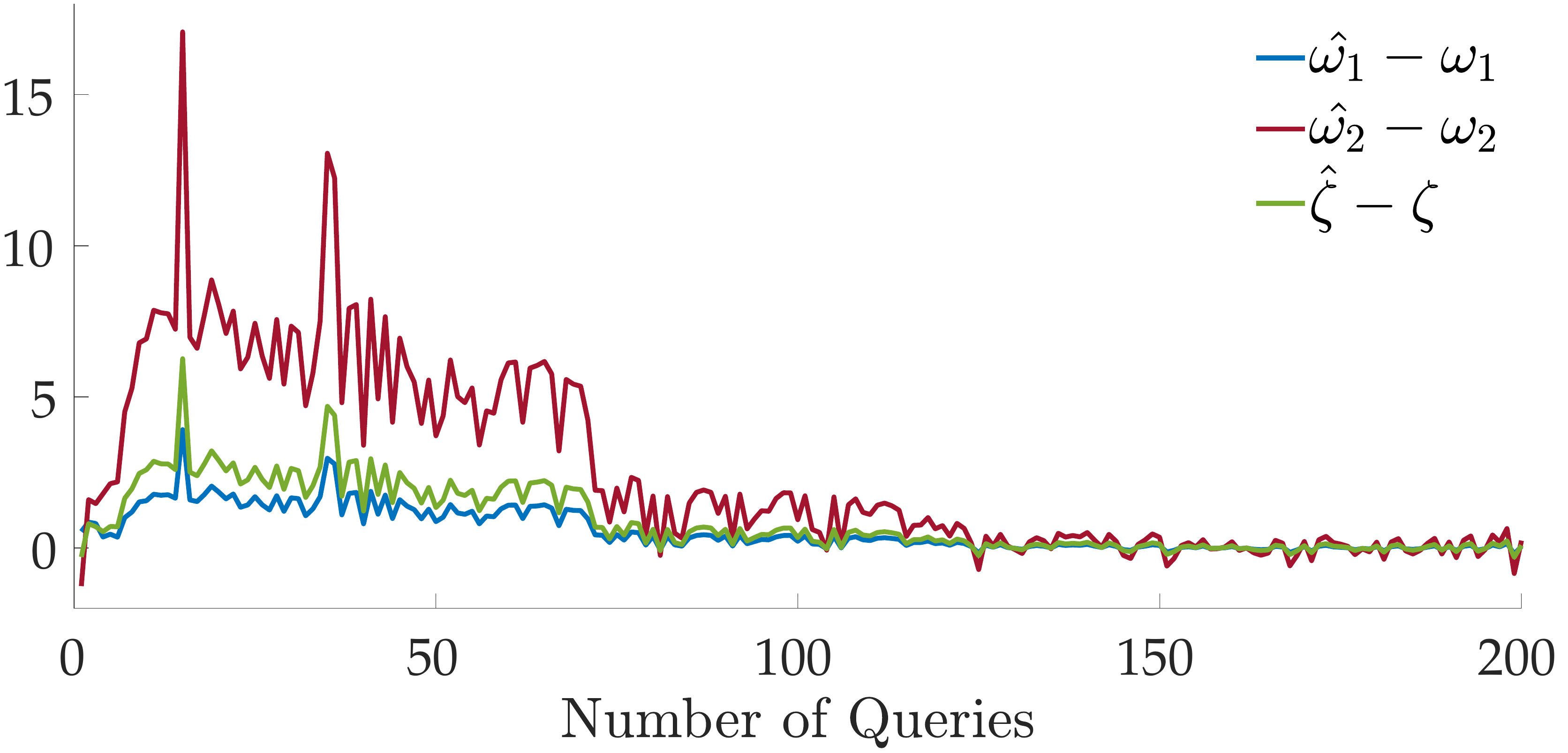}
	\caption{The errors of the reward function estimates are shown with varying number of queries. $\hat{\omega_1}$, $\hat{\omega_2}$ and $\hat{\zeta}$ represent the parameters of the sample with the highest likelihood.}
	\vspace{-5px}
	\label{fig:param_values}
\end{figure}

Fig.~\ref{fig:param_values} shows how the estimates evolved within active learning setting for one of the users. All values are overestimated initially. This is intuitively because getting noiseless responses has higher likelihood. As we query more, accepting some of the responses as noisy maximizes the likelihood. Therefore, the values start decreasing.

\begin{figure}[h]
	\centering
	\vspace{-5px}
	\includegraphics[width=\linewidth]{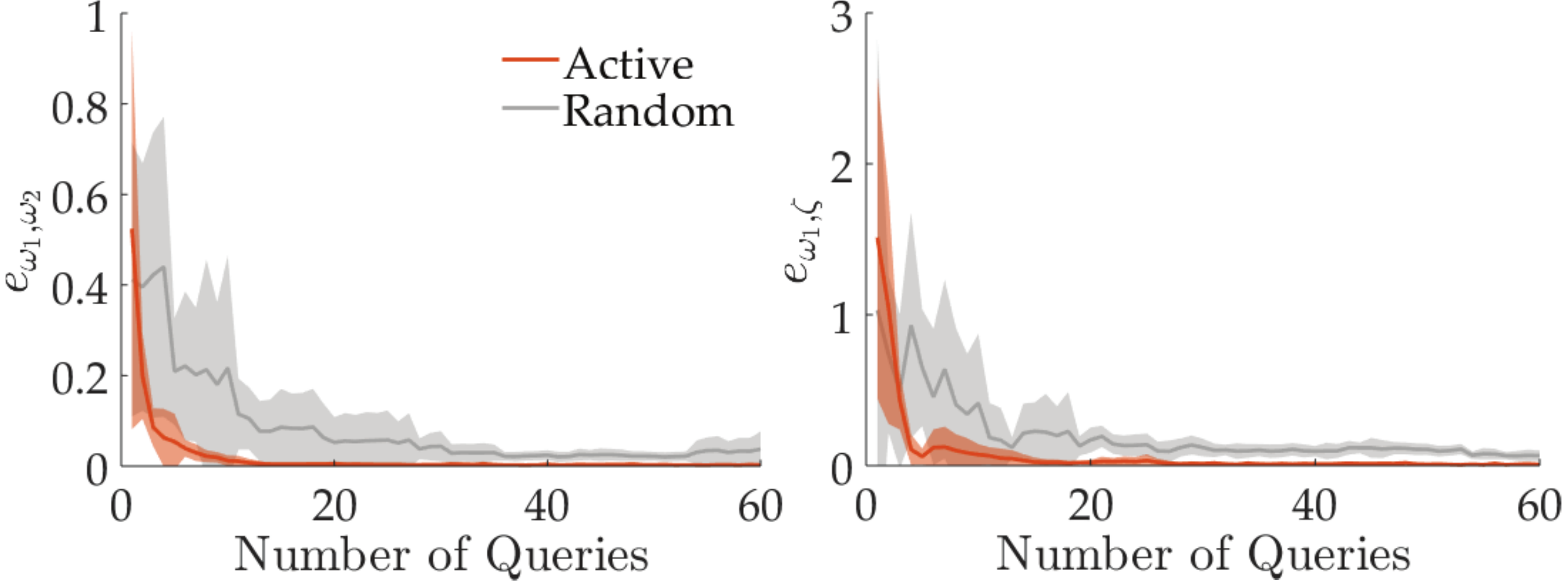}
	\caption{The error metric is averaged over $5$ different reward functions, and the results are plotted. Active learning framework provides a data-efficient way of learning by converging fast.}
	\vspace{-5px}
	\label{fig:active_vs_random}
\end{figure}

Another important observation is that the estimates of the parameters increase and decrease together even in the early iterations. This suggests we are able to learn the ratio between the parameters, e.g. $\omega_1/\omega_2$ very quickly. To check this claim, we used the following error metric:
\begin{align*}
e_{x,y} = \norm{x/y - \hat{x}/\hat{y}}_1
\end{align*}
where $x,y\in\{\omega_1,\omega_2,\zeta\}$ and $\hat{x},\hat{y}$ stand for the parameters $x, y$ of the sample with the highest likelihood. Figure~\ref{fig:active_vs_random} shows how this error decreases with increasing number of queries. It also shows how active querying enables data-efficient learning compared to the random querying baseline.

This indicates we are able to learn the relationship between parameters even under $20$ queries with active learning. All these results strongly support \textbf{H1}. 

The fact that we are able to learn the ratios implies we can estimate which road the user is most likely to choose. We will only be unsure about how noisy the user is if the parameter estimates did not converge yet. Therefore, we can still use the estimates for our planning optimization formulation even when we have small number of queries.

To test \textbf{H2}, we adopt $3$ of the equilibria from \cite{biyik2018altruistic} for comparison, namely, Nash Equilibrium where all drivers are selfish; Best Nash Equilibrium which gives the smallest average latency again with selfish drivers; and Fully Altruistic Best Nash Equilibrium which gives the smallest average latency where autonomous users are fully altruistic and human-driven cars route selfishly. Here, we give the average latencies for the specific $4$-road network from that study which we visualize in Fig.~\ref{fig:wafr_network}, and where $\demandhum=0.4, \demandaut=1.2$ cars per second:
\begin{itemize}
	\item Nash Equilibrium: $400.00$ seconds
	\item Best Nash Equilibrium: $125.66$ seconds
	\item Fully Altruistic Best Nash Equilibrium: $102.85$ seconds
\end{itemize}

We then assumed we perfectly learned the preferences of the $5$ simulated users. We ran the planning optimization with $\bar{P}=0$ and $3$ different $\theta$ to show the trade-off. The results are summarized in Table~\ref{tab:wafr_new}.

\begin{table}[H]
	\caption{Results of Routing Simulation}
	\label{tab:wafr_new}
	\centering
	\begin{tabular}{ccc}
		\hline
		$\theta$ & Avg. Latency (seconds) & Flow (cars/second) \\\hline
		$1$ & $90.41$ & $0.4412$ \\
		$20$ & $97.03$ & $1.2746$ \\
		$10^6$ & $111.28$ & $1.5964$ \\\hline
	\end{tabular}
	\vspace{-10px}
\end{table}

It can be seen we can adjust the trade-off between average latency and the served flow by tuning $\theta$. Also, given the human preferences, even when we served (almost) all of autonomous demand, our framework outperforms Best Nash Equilibrium. This shows the effectiveness of our framework on creating altruism and supports \textbf{H2}.

For \textbf{H3}, we recruited $21$ subjects ($9$ female, $12$ male) with an age range from 19 to 60. To increase the diversity among user preferences, the subjects are selected from a diverse set of locations: California / Georgia / Illinois, USA; Ankara / Denizli, Turkey; Sao Paulo, Brazil; Singapore.

In the first phase of the experiment, each participant was asked $40$ queries ($4$ roads + $1$ walking option) which are all actively synthesized. No user has selected a strictly dominated option in any query, which indicates the data are reliable. We then used their responses to get the maximum likelihood estimates which we approximate as the sample $(\boldsymbol{\bar\omega},\bar{\zeta})$ that has the highest likelihood. Afterwards, we designed $5$ different road networks each with $4$ different roads and an additional route where people may choose to walk. The roads on the $5$ different networks cover a range of different road lengths from $1.8$ kilometers to $78$ kilometers. For each network, we also set different $\theta$, $\demandhum$, $\demandaut$, and $\bar{P}$. By assuming all autonomous flow is in the service of these $21$ subjects, or of the groups that match with their preferences, we locally solved the planning optimization to obtain the pricing scheme for each traffic network. We refer to the results of this optimization as \emph{anticipated} values, e.g. anticipated profit.

In the second phase of the user study, we presented the route-price pairs and the walking option to the same $21$ subjects. For each of the $5$ networks, they picked the option they would prefer. Using these responses, we allocated the autonomous flow into the roads. However, it is technically possible that more users select a road than its maximum flow. To handle such cases, we assumed extra flow is transferred to the roads with smaller latencies without making the users pay more. If that is not feasible, extra flow is transferred to the slower roads, without any discount. While these break the fairness constraint, it rarely happens and affects only a very small portion of the flow. After autonomous flows are allocated, human driven cars selfishly chose their routes in a way to minimize the overall average latency. We refer to the results of this allocation as \emph{actual} values, e.g. actual profit.

Table~\ref{tab:user_study} compares the anticipated and the actual values. We report latencies in seconds, flows in cars per second, and profit is a rate value with the unit USD per second. In order to show how our framework incentivizes altruistic behavior, we also added two other selfish best Nash equilibrium baselines: one where the same flow as actual flow is routed (BNE1) and one where all flow demand is routed (no walking option) again under best Nash equilibrium (BNE2).

\begin{figure}[h]
	\centering
	\vspace{-5px}
	\includegraphics[width=\linewidth]{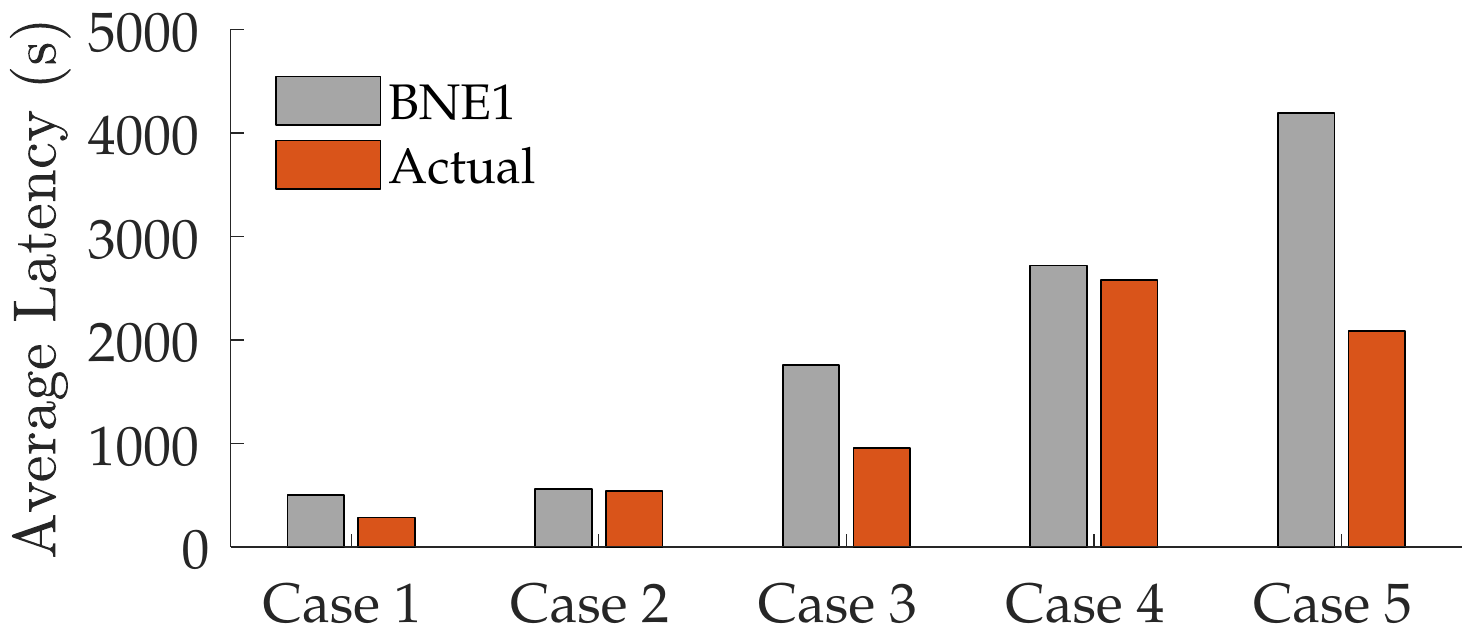}
	\vspace{-15px}
	\caption{The comparison of the actual results and BNE1, both of which allocate the same amount of flow. Our framework consistently reduces the average latency.}
	\vspace{-10px}
	\label{fig:actual_vs_bne1}
\end{figure}

It can be seen that there is generally an alignment between anticipated and actual values. While the mismatch may be further reduced by doing more queries or having more users, the difference with BNE methods is significant. In all cases, our framework achieved to incentivize altruism, which yielded lower average latencies compared to BNE1. Especially in Case 3 and Case 5, our framework approximately halved the average latency compared to the best Nash equilibria. We visualize this difference in Fig.~\ref{fig:actual_vs_bne1}. Furthermore, our framework successfully reduced flow demand when satisfying the full demand is impossible (Case 1).

One caveat is the small amount of actual flow in Case 1, which also caused an important profit loss. This is because the roads are relatively shorter, and most users preferred walking over paying for an autonomous car. Our framework could not predict this, because the learned reward functions failed to accurately estimate the probabilities. This issue will be further discussed in Section~\ref{sec:conclusion}.

\begin{table*}[h]
	\caption{Results of Real-User Experiments}
	\label{tab:user_study}
	\centering
	\begin{tabular}{lccccccccccc}
		\hline
		& \multicolumn{3}{c}{\textbf{Anticipated}} & \multicolumn{3}{c}{\textbf{Actual}} & \multicolumn{2}{c}{\textbf{BNE1}} & \multicolumn{2}{c}{\textbf{BNE2}}  \\ \hline
		& Avg. Latency & Flow & Profit & Avg. Latency & Flow & Profit & Avg. Latency & Flow & Avg. Latency & Flow \\ \hline 
		Case 1 & $338.10$ & $1.6622$ & $11.11$ & $283.30$ & $1.0429$ & $3.85$ & $501.97$ & $1.0429$ & \multicolumn{2}{c}{Infeasible} \\
		Case 2 & $510.97$ & $1.5592$ & $4.50$ & $537.30$ & $1.5048$ & $3.95$ & $562.20$ & $1.5048$ & $562.20$ & $1.6000$ \\
		Case 3 & $921.06$ & $1.2677$ & $4.59$ & $957.71$ & $1.3000$ & $3.99$ & $1756.89$ & $1.3000$ & $1756.89$ & $1.3000$ \\
		Case 4 & $2083.42$ & $1.2383$ & $48.48$ & $2089.53$ & $1.2048$ & $46.32$ & $4194.26$ & $1.2048$ & $4194.26$ & $1.3000$ \\
		Case 5 & $2576.36$ & $1.5538$ & $40.01$ & $2579.52$ & $1.6000$ & $49.33$ & $2720.00$ & $1.6000$ & $2720.00$ & $1.6000$ \\ \hline
	\end{tabular}
	\vspace{-15px}
\end{table*}

\section{Conclusion}
\label{sec:conclusion}
In this work, we develop a method improving the efficiency of traffic networks shared between human-driven and autonomous vehicles by influencing the routing behavior of humans. Considering an autonomous service that offers a choice of a number of routes, we formulate an optimization to find the prices that influence autonomous user choices so that when users choose their routes, and the human-driven traffic routes solely based on travel delay, travel delay is minimized and road usage is maximized. In order to do so, we model how people choose between between routes with differing prices and latencies, and develop a method to actively learn the choice model parameters for a population of users. We verify these mechanisms through an experiment with real human choice data.

There are areas in which our mechanism can be improved. Specifically, one could further refine the reward function model for autonomous service users and not necessarily assume that it's linear. One could also consider the human-driven vehicle demand to be elastic, meaning that people may choose to stay home or choose public transportation if congestion is high. Further, we find that when solving the planning optimization using actual data from user choices, the profit constraint is often unsatisfied. To manage this, we could instead include profit as part of the objective function or optimize for the worst-case performance.

We considered only a single alternative to the autonomous service use. However, a population may choose from a variety of options such as walking, biking, or taking the bus, $\zeta$ will have a multimodal distribution; we then need to learn the mixture. Finally, it is unclear if it is better to use many samples from our posterior on $(\boldsymbol{\omega},\zeta)$ or just the maximum likelihood sample while computing $\mathbf{q}$ for the planning optimization.

Future work can address these issues or expand in broader directions. It would be worthwhile to look at more general network topologies. Further, the information that people have access plays a huge role in decision making. It would be enlightening to develop choice models that incorporate this information availability as well as possibly other aspects of decision making, such as risk aversion. These expansions would do much to address congestion on roads that are soon to be shared between human drivers and autonomous vehicles.

\bibliographystyle{IEEEtran}
\bibliography{IEEEabrv,paper}

\section*{Appendix}
\subsection{Proof of Proposition \ref{prop:eq_road}}
Let us call the reward function that is homogeneous among the users as $\bar{r}$. Assume the optimal solution is $(\latencyvec^*, \pricevec^*, {\flowhumvec}^*, {\flowautvec}^*)$, and $\latencyvec^*_k > a_k$, where $k=\argmax_\pathidx \fflatency_\pathidx$ subject to ${\flowhum}^*_{\pathidx} > 0$. Due to selfishness, we know $\latency_{\pathidx}=\latency_k$ for $\forall \pathidx\in[k]$. Also, let $D^*$ (and $\bar D^*$) be the set of dominated (undominated) options under the optimal solution.

Let $(\latencyvec',\pricevec')$ be such that
\begin{itemize}[nosep]
	\item $\latency'_{\pathidx}=\latency^*_{\pathidx}$ and $\price'_{\pathidx}=\price^*_{\pathidx}$, for $\forall \pathidx\in\pathset\setminus[k]$,
	\item $\fflatency_k \leq \latency'_k < \latency^*_k$,
	\item $\latency'_\pathidx = \latency'_k$, for $\forall \pathidx\in[k-1]$,
	\item $\price'_\pathidx = \price^*_\pathidx + \epsilon$ where $\epsilon>0$ for $\forall \pathidx \in [k]\cup D^*$,
	\item $\mathbf{q}(\latencyvec',\pricevec') = \mathbf{q}(\latencyvec^*,\pricevec^*)$ for $\pathidx\in[k]\; .$
\end{itemize}
We also define $D'$ and $\bar{D}'$ similarly. The only two free variables of $(\latencyvec',\pricevec')$ are $\latency'_k$ and $\epsilon$. To show such an $(\latencyvec',\pricevec')$ with the last constraint exists, we prove
1) For any $\latency'_k$, there exists an $\epsilon$ such that $\bar{r}(\latencyvec^*,\pricevec^*,\pathidx)=\bar{r}(\latencyvec',\pricevec',\pathidx)$ for $\forall\pathidx\in[k]$, 2) $D^*=D'$.
These together imply $\bar{r}(\latencyvec^*,\pricevec^*,\pathidx)=\bar{r}(\latencyvec',\pricevec',\pathidx)$ for $\forall \pathidx\in\pathset$, which then implies the last constraint.

The first one is easy to show. $\bar{r}(\latencyvec^*,\pricevec^*,\pathidx)=\bar{r}(\latencyvec^*,\pricevec^*,k)$ for $\forall \pathidx\in[k]$. For any $\latency'_k$, we can set $\epsilon=\frac{\omega_1}{\omega_2}(\latency^*_k-\latency'_k)$ to get $\bar{r}(\latencyvec',\pricevec',i)=\bar{r}(\latencyvec^*,\pricevec^*,\pathidx)$ for $\forall \pathidx \in[k]$ where $\boldsymbol{\omega}$ is associated with $\bar{r}$. The equality also holds for $\pathidx\in\pathset\setminus[k]$, because the prices and latencies are the same between $(\latencyvec^*,\pricevec^*)$ and $(\latencyvec',\pricevec')$ in those roads.

We show $D^*=D'$ in two steps:\\
1) $D^*\cap[k] = D'\cap[k]$, since the roads had the same latency drop and price increase, and they cannot be dominated by any of $\pathset\setminus[k]$, because they have strictly lower latency.\\
2) $D^*\setminus[k] = D'\setminus[k]$, because the roads have the same latency between two configurations and the prices went up only for those in $D^*$, so they cannot dominate each other in a different way between the configurations. They also cannot dominate any of the roads in $[k]$, because they have strictly higher latency. And $\pathidx\in D'\setminus[k]\implies\pathidx\in D^*\setminus[k]$, because the roads in $[k]$ had already lower latency in the optimal solution and their prices only went up.

Having showed $(\latencyvec',\pricevec')$ exists, we note $(\latencyvec',\pricevec',{\flowhumvec}',{\flowautvec}^*)$ is a feasible solution for some ${\flowhumvec}'$, because the distribution of autonomous service users remained the same and the prices went only up, which together means the profit constraint \eqref{c:7} is still satisfied. The other constraints are trivially satisfied by the construction of $(\latencyvec',\pricevec')$ and as the optimal solution has to be feasible. Then, since the demand of human-driven cars are inelastic, the total served flow is the same between two solutions. However, the latencies of the roads in $[k]$ decreased, which means their capacities increased. Hence, $k\!\geq\!\argmax_\pathidx a_i$ subject to ${\flowhum}'_{\pathidx}\!>\!0$. Then, the average latencies of human-driven cars decreased. Similarly, the average latencies of autonomous service users may only go down. Hence, $J({\flowhumvec}^*,{\flowautvec}^*,\congprofvec^*)\!>\! J({\flowhumvec}',{\flowautvec}',\congprofvec')$, and $(\latencyvec^*, \pricevec^*, {\flowhumvec}^*, {\flowautvec}^*)$ cannot be the optimal solution.$\qquad\square$

\end{document}